\documentstyle[12pt]{article}
\newtheorem{defi}{Definition}

\def\today{\ifcase\month\or
January\or February\or March\or April\or May\or June\or July\or August\or September\or October\or November\or December\fi
\space\number\day ,\number\year}

\newcommand{\proof}{\noindent{\bf Proof}~~}
\newtheorem{Theorem}{Theorem}
\newtheorem{Proposition}{Proposition}
\newtheorem{Lemma}{Lemma}
\newcommand{\bl}{\begin{Lemma}}
\newtheorem{Remark}{Remark}
\newcommand{\el}{\end{Lemma}}
\newcommand{\be}{\begin{equation}}
\newcommand{\ee}{\end{equation}}
\newcommand{\bd}{\begin{defi}}
\newcommand{\ed}{\end{defi}}

\newtheorem{pro}{Proposition}
\newcommand{\bp}{\begin{pro}}
\newcommand{\ep}{\end{pro}}
\newcommand{\bt}{\begin{Theorem}}
\newcommand{\et}{\end{Theorem}}
\newtheorem{cor}{Corollary}
\newcommand{\bc}{\begin{cor}}
\newcommand{\ec}{\end{cor}}

\def\sqr#1#2{{\vcenter{\vbox{\hrule height.#2pt
\hbox{\vrule width.#2pt height#1pt \kern#1pt
\vrule width.#2pt}\hrule height.#2pt}}}}
\def\square{\mathchoice\sqr45\sqr45\sqr{2.1}3\sqr{1.5}3}
\newcommand{\qed}{\square}

 at 10truept

\begin{document}
\setlength{\textheight}{7.7truein}    
\setcounter{page}{1} \centerline{\bf Generalized Eta-Einstein and $(\kappa ,\mu )$-structures} \centerline{\bf } \baselineskip=13pt
\vspace*{10pt} \centerline{\footnotesize {\bf 
Philippe Rukimbira}} \baselineskip=12pt \centerline{\footnotesize\it
Department of Mathematics and Statistics, Florida International University}
\baselineskip=10pt \centerline{\footnotesize\it Miami, Florida
33199, USA} \baselineskip=10pt \centerline{\footnotesize E-MAIL:
rukim@fiu.edu} \vspace*{0.225truein}

\vspace*{0.21truein} \abstract{\it   Geneneralized $(\kappa ,\mu )$- structures occur in dimension 3 only. In this dimension 3,  only K-contact structures can occur as  generalized Eta-Einstein.
On closed manifolds, Eta-Einstein, K-contact structures which are not D-homothetic to Einstein structures  are almost regular. We also construct examples of compact, generalized Jacobi $(\kappa ,\mu )$-structures.} {}


\vspace*{14pt}                  

\baselineskip=24pt


	
\section{Introduction} 
 The purpose of this paper is to study the interaction between generalized $(\kappa ,\mu )$-structures and generalized Eta-Einstein structures. In 1995, Blair, Koufogiorgos and Papantoniou \cite{BKP}  introduced and studied extensively $(\kappa ,\mu )$-spaces. Generalized $(\kappa ,\mu)$ where at least one of $\kappa$ and $\mu$ is  assumed not to be a constant function were studied in many papers, including \cite{KOC} where it is shown that there are no generalizations in odd dimensions 5 or higher. We use the terminology " generalized" when at least one of the functions $\kappa $ and $\mu$ is not constant.

Generalized Eta-Einstein structures were introduced and studied by Okumura \cite{OKU} in 1962. An extensive study of these structures from the Sasakian point of view has been carried out  by Boyer, Galicki and Matzeu \cite{BOG}. As it turns out, in dimension 3, Eta-Einstein structures are exactly $(\kappa ,\mu )$-structures with $\mu =0$, simply referred to as $\kappa$-structures. One of the results in this paper is that non-K-contact generalized Eta-Einstein structures do not occour in dimension 3.

Using a new type of contact metric deformation, we construct examples of compact, generalized Jacobi $(\kappa ,\mu )$ structures.

\section{Basic notions of contact metric structures}

A contact form on a $2n+1$-dimensional manifold $M$ is a one-form $\eta$ such that $\eta\wedge (d\alpha )^n$ is a volume form on $M$. Given a contact manifold $(M,\eta )$, there exist tensor fields $(\xi ,\phi , g)$, where $g$ is a Riemannian metric and $\xi$ is a unit vector field, called the Reeb field of $\eta$ and $\phi$ is an endomorphism of the tangent bundle of $M$ such that
\begin{itemize}
\item[(i)] $\eta (\xi )=1,~\phi^2 =-Id +\eta\otimes\xi ,~\phi\xi =0$
\item[(ii)] $d\eta =2g(.,\phi .)$
\end{itemize}
The data $(M,\eta ,\xi ,\phi ,g)$ is called a contact metric structure on $M$.

Denoting by $\nabla$ the Levi-Civita connection of $g$, and by $$R(X,Y)Z=\nabla_X\nabla_YZ-\nabla_Y\nabla_XZ-\nabla_{[X,Y]}Z$$ its curvature tensor, a contact metric structure $(M,\eta ,\xi ,\phi ,g)$ is called Sasakian if the condition$$(\nabla_X\phi )Y=g(X,Y)\xi -\eta (Y)X$$ is satisfied for all tangent vectors  $X$ and $Y$.
A well known curvature characterization of the Sasakian condition is as follows:
\begin{pro}
{\it 
A contact metric structure $(M,\eta , \xi ,\phi ,g)$ is Sasakian if and only if $$R(X,Y)\xi =\eta (Y)X-\eta (X)Y$$ for all tangent vectors $X$ and $Y$.}\end{pro}

Contact metric structures on which the Reeb vector field is an infinitesimal isometry are called K-contact structures. One has the following characterization of K-contact structures.
\begin{pro}{\it A contact metric structure $(M,\eta ,\xi ,\phi ,g)$ is K-contact if and only if its curvature tensor$R$ satisfies $$R(X,\xi )\xi =X-\eta (X)\xi$$ for all tangent vectors $X$.}\end{pro}

For more on the contact metric geometry, reference \cite{BLA} includes an extensive bibliography and constitute a very good introduction to the subject. 

We will need the following fact related to contact distributions being bracket generating. First, some terminology being used throughout the paper. A tangent vector $U$ on a contact manifold $(M,\eta ,\xi  )$ is said to be horizontal is $\eta (U)=0$.
\begin{Lemma} \label{lem1m} If a smooth function $f$ on a contact manifold $(M,\eta ,\xi )$ satisfies $df(U)=0$ for every horizontal vector $U$; then $df=0$ identically and $f$ is locally constant.
\end{Lemma}

\begin{proof}  The condition $df(U)=0$ for all  horizontal $U$  implies that $df=df(\xi)\eta$. We claim that $df(\xi)$ is identically zero. If not, the contact distribution would be then in the kernel of the  closed form $df$, hence integrable, contradicting the definition of contact distributions as totally non-integrable. Therefore, $df=df(\xi)\eta =0$.
\end{proof}

\section{Generalized ($\kappa , \mu$ )-structures}

A contact  metric space $(M,\eta ,\xi ,\phi .g)$ is said to be  generalized $(\kappa ,\mu )$ if its (1,3)-curvature tensor $R$ satisfies$$R(X,Y)\xi =\kappa (\eta (Y)X-\eta (X)Y)+\mu (\eta(Y)hX-\eta(X)hY)$$ for some functions $\kappa$ and $\mu$ and the tensor field $h$ is difined by $h=\frac{1}{2}L_\xi\phi$. By generalized $(\kappa ,\mu )$, we mean that at least one of the functions $\kappa$ and $\mu$ is not constant. In dimensions $2n+1>3$, the functions $\kappa$ and $\mu$ are known to be necessarily constant \cite{KOC}. In dimension 3, non-compact examples of   generalized $(\kappa , \mu )$ structures have been presented in \cite{KOC}. Moreover, in \cite{SHA}, Sharma has shown that on $\kappa$-spaces, i.e $(\kappa ,\mu )$ with $\mu =0$; the function $\kappa$ is necessarily constant. Here, we prove a slight generalization of this result under Proposition \ref{prop1m} further in the paper.

In general, on any Riemannian 3-dimensional manifold $(M,g)$, denoting by $Q$, $r$ and $Ric (,)$ the Ricci operator, the scalar curvature and the Ricci tensor respectively, the following identity holds \cite{PET}: 

  $$ R(X,Y)Z=$$ \begin{equation}Ric(Z,Y)X-Ric(Z,X)Y+g(Z,Y)QX-g(Z,X)QY\\ -\frac{r}{2}\{g(Z,Y)X-g(Z,X)Y\} \label{equatpet}\end{equation}
for any vector fields $X$, $Y$, and $Z$ on $M$.

On a 3-dimensional $(\kappa ,\mu )$ space $(M,\eta , \xi , \phi ,g)$, letting $Z=\xi =Y$  in identity (\ref{equatpet}) above,
we obtain the following:
\begin{equation}\label{equatm2}R(X,\xi )\xi =Ric (\xi ,\xi )X-Ric (\xi , X)\xi +QX-\eta (X)Q\xi -\frac{r}{2}(X-\eta (X)\xi )\end{equation} 

 Using the fact that on a $(\kappa ,\mu )$-space, $Q\xi =2\kappa\xi$ in (\ref{equatm2}) and solving for $QX$,  we obtain:
\begin{equation}\label{equat2w}QX=(\frac{r}{2}-\kappa )X+(3\kappa -\frac{r}{2})\eta (X)\xi +\mu hX\end{equation}
Similar to the standard $(\kappa ,\mu )$-structures case, the eigenfunctions  of the symmetric tensor $h$ are $\lambda =\sqrt{1-\kappa}$ and $-\lambda$ on a generalized $(\kappa ,\mu )$-space with $\kappa <1$ (See \cite{KOC}). Denoting by $E$ a unit eigenvector field corresponding to $\lambda$, we have the following lemma:

\begin{Lemma}\label{keylem} On every 3-dimensional generalized $(\kappa ,\mu )$ structure $(M,\eta ,\xi ,\phi ,g)$, the following identities hold:
\begin{equation}d\kappa (\xi )=0,~~~~~dr(\xi )=0\label{ll33}\end{equation}
\begin{equation}\nabla_EE=\frac{d\lambda (\phi E)}{2\lambda}\phi E,~~~~~~\nabla_{\phi E}\phi E=\frac{d\lambda (E)}{2\lambda}E\label{ll4}\end{equation}
\begin{equation}d\mu (E)=-2d\lambda (E),~~~~~~d\mu (\phi E)=2d\lambda (\phi E)\label{ll5}\end{equation}
\end{Lemma}

\begin{proof} On any generalized $(\kappa ,\mu )$-space, the following identities hold:
\begin{equation}\nabla_\xi h=\mu h\phi ,~~~~~~h^2=(\kappa -1)\phi^2\label{ll1}\end{equation}

Differentiating the second equation along $\xi$, then substituting in the first,  one gets:
\begin{equation}\mu h^2\phi -\mu h^2\phi =d\kappa (\xi )\phi^2\end{equation} Hence $d\kappa (\xi )\phi^2=0$ which implies  $d\kappa (\xi )=0$, proving first of  (\ref{ll33}).

By $Ric$ we denote the Ricci tensor, $Ric (X,Y)=g(QX,Y)$ for any tangent vectors $X$ and $Y$.
Using identities 
\begin{equation} div~Ric=\frac{1}{2}dr\label{key3}\end{equation} valid on any Riemannian manifold and, for any frame $\{E_1, E_2, E_3\}$,
\begin{equation}\sum_{i=1}^3(\nabla_{E_i}h)E_i=\phi Q \xi\end{equation} valid on any contact metric 3-dimensional manifold (see \cite{KOU}), we obtain
\begin{equation}(\nabla_Eh)E+(\nabla_{\phi E}h)\phi E+(\nabla_\xi h)\xi =\phi Q\xi\end{equation}
On a generalized $(\kappa ,\mu )$-space, the above reduces to \begin{equation}(\nabla_Eh)E+(\nabla_{\phi E}h)\phi E=0\end{equation}Which is simply:
\begin{equation}\nabla_E(hE)-h\nabla_EE+\nabla_{\phi E}(h\phi E)-h\nabla_{\phi E}\phi E=0\label{ll2}\end{equation}

Next, using $hE=\lambda E$ and $h\phi E=-\lambda \phi E$, (\ref{ll2}) becomes:
\begin{equation}d\lambda (E)E-2\lambda\nabla_{\phi E}\phi E-d\lambda (\phi E)\phi E+2\lambda\nabla_EE=0\label{ll3}\end{equation}
Inner product (\ref{ll3}) with $ E$ shows that $\nabla_{\phi E}\phi E=\frac{d\lambda (E)}{2\lambda}E$ and inner product of (\ref{ll3}) with $ \phi E$ shows that $\nabla_EE=\frac{d\lambda (\phi E)}{2\lambda}\phi E$, proving (\ref{ll4}).

To prove (\ref{ll5}), we use the well known Riemannian geometric identity (\ref{key3}). Denoting by $grad~r$ the gradient vector of function $r$ and using the frame $\{E,\phi E,\xi\}$, identity (\ref{key3}) becomes:
\begin{equation}\frac{1}{2}grad~r=(\nabla_EQ)E+(\nabla_{\phi E}Q)\phi E+(\nabla_\xi Q)\xi \end{equation} or equivalently
\begin{equation}\frac{1}{2}grad~r=\nabla_EQE-Q\nabla_EE+\nabla_{\phi E}Q\phi E-Q\nabla_{\phi E}\phi +\nabla_\xi Q\xi-Q\nabla_\xi \xi\label{ll11}\end{equation}

From (\ref{equatm2}), we have :
\begin{equation}QE=(\frac{r}{2}-\kappa )E+\mu\lambda E\label{ll6}\end{equation}
\begin{equation}Q\phi E=(\frac{r}{2}-\kappa )\phi E-\mu\lambda \phi E\label{ll7}\end{equation}
\begin{equation}Q\xi =2\kappa \xi\label{ll8}\end{equation}

Substituting these into (\ref{ll11})and using (\ref{ll4}), the right hand side of (\ref{ll11}) becomes:
\begin{equation}\label{key1}
\begin{array}{lr}\frac{1}{2}(dr(E)E+dr(\phi E)\phi E)-d\kappa (E)E-d\kappa (\phi E)\phi E+(\frac{r}{2}-\kappa )(\frac{d\lambda(\phi E)}{2\lambda})\phi E+&\\ 
d(\mu\lambda )(E)E+\mu\lambda (\frac{\lambda (\phi E)}{2\lambda}\phi E-\frac{d\lambda(\phi E)}{2\lambda}((\frac{r}{2}-\kappa )\phi E-\mu\lambda \phi E)+(\frac{r}{2}-\kappa )(\frac{d\lambda (E)}{2\lambda}) &\\ -d(\mu\lambda )(\phi E)\phi E-(\mu\lambda )(\frac{d\lambda (E)}{2\lambda} )E-\frac{d\lambda (E)}{2\lambda}((\frac{r}{2}-\kappa )E+\mu\lambda E))\end{array}\end{equation}
while the left hand side is:
\begin{equation}\label{key2}\frac{1}{2}(dr(E)E+dr(\phi E)\phi E+dr(\xi )\xi)\end{equation}
Cancelling identical terms from(\ref{key1}) and (\ref{key2}), one reaches
\begin{equation}dr(\xi )=0\label{key12}\end{equation}
\begin{equation}\label{key14}-d\kappa (E)+d(\mu\lambda )(E)-\mu\lambda\frac{d\lambda (E)}{\lambda}=0\end{equation}
\begin{equation}\label{key15}-d\kappa (\phi E)+\mu\lambda\frac{d\lambda (\phi E)}{\lambda}-d(\mu\lambda )(\phi E)=0\end{equation}
Identity (\ref{key12}) proves the second identity of (\ref{ll33}).

Since $\kappa =1-\lambda^2$, identity (\ref{key14}) becomes $0=2\lambda d\lambda (E)+\mu d\lambda (E)+\lambda d\mu (E)-\mu d\lambda (E)=2\lambda d\lambda (E)+\lambda d\mu(E)$. Hence, since $\lambda \ne 0$, one gets $
d\mu (E)=-2d\lambda (E)$, proving the first identity of (\ref{ll5}). In the same way, identity (\ref{key15}) becomes $d\mu (\phi E)=2d\lambda (\phi E)$, proving the second identity of (\ref{key15}).
\end{proof}
$\qed$
 
We may now prove the following:
\begin{pro}\label{prop1m} Let $(M,\eta , \xi ,\phi ,g)$ be a 3-dimensional generalized $(\kappa ,\mu )$-space with  $$R(X,Y)\xi =\kappa (\eta (Y)X-\eta (X)Y)+\mu (\eta(Y)hX-\eta(X)hY$$ for any tangent vectors $X,Y$  with $\kappa <1$.  Then the function $\kappa$ is constant if and only if the function $\mu$ is constant. In particular, on any $\kappa$-space,  $\mu =0$, and hence $\kappa$ is constant.
\end{pro}

\begin{proof}

If $\mu$ is constant, then $d\lambda (E)=0=d\lambda (\phi E)$ from (\ref{ll5}) in Lemma \ref{keylem}. Hence by Lemma \ref{lem1m}, $d\lambda =0$ and $\lambda$ is constant. But since $\lambda  =\sqrt{1-\kappa}$, it follows that $\kappa$ is also constant. Conversely, if $\kappa$ is constant,then $\lambda$ is constant and a similar argument shows that $\mu$ is also constant.  $\qed$
\end{proof}

\section{Generalized Eta-Einstein structures}

A contact metric space $(M,\eta ,\xi ,\phi ,g)$ is said to be generalized Eta-Einstein if its Ricci operator $Q$ satisfies $$QX=\lambda X+\gamma \eta (X)\xi$$ for some smooth functions $\lambda $ and $\gamma$. In dimensions $2n+1>3$, it is known that when $Ric (\xi, \xi )=g(Q\xi ,\xi )$ is constant, like in the case of K-contact generalized $\eta$-Einstein, then $\lambda$ and $\gamma$ are constant functions.

More generally, we prove the following:
\begin{pro} \label{prop2m}Let $(M,\eta ,\xi ,g)$ be a generalized Eta-Einstein $2n+1$-dimensional manifold.Then the Ricci operator $Q$ is given by $$QX= \lambda X+\gamma \eta(X)\xi$$ where $d\lambda(\xi )=0=d\gamma (\xi )$ and $$(1-2n)d\lambda =d\gamma$$ Moreover any of the functions $\lambda$ and $\gamma$ is constant if and only if the other one is.
\end{pro}

\begin{proof}   Denoting the scalar curvature by $r$ and the Ricci curvature tensor  by $Ric$, then one has $$r=(2n+1)\lambda +\gamma$$ and $$div Ric(X)=d\lambda(X) +d\gamma(\xi)\eta(X)$$

From the well known Riemannian geometric identity$$2 div Ric=dr$$ we get
$$2d\lambda (X)+2d\gamma(\xi)\eta(X)=(2n+1)d\lambda (X)+d\gamma(X)$$
Equivalently
\begin{equation}\label{equat9w}(1-2n)d\lambda(X)=d\gamma(X)-2d\gamma (\xi)\eta(X)\end{equation}
Evaluating (\ref{equat9w}) on  $H$ orthogonal to $\xi$, we get \begin{equation}(1-2n)d\lambda (H)=d\gamma (H)\label{equat2m}\end{equation} Identity (\ref{equat2m}) implies that \begin{equation}\label{equat10w}(1-2n)d\lambda =d\gamma \end{equation}  by Lemma \ref{lem1m}.

Letting $X=\xi$ in  (\ref{equat9w}), we get
$$(1-2n)d\lambda(\xi)=-d\gamma(\xi)$$ This combined with (\ref{equat10w}) leads to $$d\lambda (\xi)=0=d\gamma (\xi)$$
$\qed$
\end{proof}

In the 3-dimensional setting, we have the following theorem:
\begin{Theorem} \label{theorem1w}Let $(M,\eta,  \xi ,g)$ be 3-dimensional, generalized Eta-Einstein, not K-contact manifold. Then the Ricci operator $Q$ is given by $$QX=\gamma \eta(X)\xi$$ with constant $\gamma$. In particular, $(M,g)$ has constant scalar curvature $\gamma$.
\end{Theorem}

\begin{proof} Suppose \begin{equation}\label{eqricci1}QX=\lambda X+\gamma\eta(X)\xi\end{equation} We will show that $\lambda =0$ and $d\gamma =0.$

The curvature tensor of the 3-dimensional, generalized  Eta-Einstein manifold satisfies \cite{PET}:
$$R(X,Y)\xi=Ric(Y,\xi)X-Ric(X,\xi)Y+\eta(Y)QX-\eta(X)QY-\frac{r}{2}(\eta(Y)X-\eta(X)Y)$$
\begin{equation}=(2\lambda+\gamma -\frac{r}{2})(\eta(Y)X-\eta(X)Y)\label{12w}\end{equation}
Identity (\ref{12w}) means that $(M,\eta,\xi,g)$ is a $\kappa$-space, with $\kappa =2\lambda +\gamma -\frac{r}{2}$; which is constant by \cite{SHA}. 

By a result in \cite{BKP}, the scalar curvature of a 3-dimensional $\kappa$-space is $r=2\kappa$. Thus, in the case at hand, one has: $$r=2(2\lambda+\gamma -\frac{r}{2});$$ That is : $r=2\lambda +\gamma$. This compared with the calculated value from (\ref{eqricci1}), $r=3\lambda +\gamma$,  leads to $\lambda =0$.

Since by Proposition \ref{prop2m},  $d\lambda$ is proportional to $d\gamma$, we see that $d\gamma =0$.

$\qed$

\end{proof}

One has the following immediate corollary of Theorem \ref{theorem1w}:

\begin{cor} If a 3-dimensional contact metric structure $(M,\eta ,\xi ,\phi ,g) $ is strict generalized Eta-Einstein, then  $(M,\eta ,\xi ,\phi ,g)$ must be K-contact.\end{cor}

\section{Automorphisms of generalized Eta-Einstein and $(\kappa ,\mu )$-structures}

The automorphism groups of generalized Eta-Einstein and $(\kappa ,\mu )$ structures exhibit some remarkable properties.  For generalized Eta-Einstein structures, we have the following:
\begin{Theorem}\label{theo1} Let $Z$ be a Killing vector field on a generalized Eta-Einstein space $(M,\eta ,\xi ,\phi ,g)$ with Ricci operator $Q=\lambda Id +\gamma \eta \otimes \xi$, $\gamma\ne 0$ almost everywhere.
Then $Z$ is an infinitesimal automorphism of the contact metric structure. Moreover, $\lambda$ and $\gamma$ are integral invariants of $Z$, that is, $d\lambda (Z)=0=d\gamma (Z)$.
\end{Theorem}

First a needed  proposition before the proof of Theorem \ref{theo1}:

\begin{pro}\label{pro1} On any contact metric structure $(M,\eta ,\xi ,\phi, g)$, if $Z$ is an infinitesimal isometry, then $\eta ([Z,\xi ])=0.$\end{pro}

\begin{proof}  $$\begin{array}{rcl}\eta ([Z,\xi ])&=&g(\xi ,[Z,\xi ])=g(\xi ,\nabla_Z\xi -\nabla_\xi Z)\\&=&\frac{1}{2}Zg(\xi ,\xi )-g(\xi ,\nabla_\xi Z)=0-0.\end{array}$$ We used the fact that $\nabla Z$ is a skew-symetric operator.

\end{proof}

\begin{proof}[Proof of Theorem~{\upshape\ref{theo1}}]

 Lie differentiating $Q=\lambda Id+\gamma\eta \otimes \xi$ in the direction of $Z$, we obtain:
\begin{equation}\label{1}0=Z(\lambda )Id +Z(\gamma )\eta \otimes \xi +\gamma ( L_Z\eta )\otimes \xi +\gamma\eta \otimes [Z, \xi ]\end{equation}
 Evaluating (\ref{1}) at $\xi$ and inner product (\ref{1}) with an arbitrary horizontal vector $V$;
one obtains:
$$\gamma g([Z,\xi ], V)=0$$ 
That is, since $[Z,\xi ]$ is horizontal by Proposition \ref{pro1}, $[Z,\xi ]=0$ and hence $L_Z\eta =0$ whenever $\gamma\ne0$. Taking this into account,  identity (\ref{1}) reduces to 
\begin{equation}\label{2R}0=Z(\lambda)Id+Z(\gamma)\eta\otimes\xi\end{equation}

Evaluating (\ref{2R}) at an arbitrary horizontal $V$ leads to $0=Z(\lambda)V$, hence $Z(\lambda)=0$. Then, evaluating  (\ref{2R}) at $\xi$ gives $Z(\gamma)=0$.

$\qed$

\end{proof}

In dimension 3, any $\kappa$-structure is  generalized Eta-Einstein and conversely, any generalized Eta-Einstein space is a $\kappa$-space. This is easily seen from identities (\ref{equat2w}) and (\ref{12w}).  In \cite{TAN}, a theorem of Tanno,  Theorem A, states that on a 3-dimensional K-contact manifold which is not of constant curvature, every infinitesimal isometry is an automorphism of the contact metric structure. We generalize this result in two directions: First to all odd dimensions, then  to 3-dimensional $\kappa$-structures.  

\begin{Theorem} Let $(M,\eta ,\xi ,\phi ,g)$ be a K-contact, generalized Eta-Einstein space which is not Einstein. Then any infinitesimal isometry is an infinitesimal automorphism.
\end{Theorem}

\begin{proof} This is a direct consequence of Theorem \ref{theo1}. The non-Einstein condition implies that the parameter $\gamma$ in the Ricci operator identity is not zero.\end{proof}

In the particular dimension 3, we have:
\begin{Theorem} Let $(M,\eta ,\xi ,\phi ,g)$ be a 3-dimensional $\kappa$-space. If $\kappa \ne 0,1$ or $\kappa =1$ and the scalar curvature $r\ne 6$, then every infinitesimal isometry is an automorphism of the $\kappa$-structure.\end{Theorem}

\begin{proof} From Identity (\ref{equat2w}), it follows that the Ricci operator $Q$  of any $\kappa$-space is given by
\begin{equation}\label{eqric}QX=(\frac{r}{2}-k)X+(3\kappa -\frac{r}{2})\eta(X)\xi\end{equation}

In the case $\kappa \ne 0,1$; we know $r=2\kappa$ and $Q$ is given by $QX=2\kappa\eta(X)\xi$ with $2\kappa \ne 0$, see \cite{BKP}.

In the case $\kappa =1$ but $r\ne 6$, it holds $QX=(\frac{r}{2}-1)X+(3-\frac{r}{2})\eta(X)\xi$ with $3-\frac{r}{2}\ne 0$.

In each case, the conclusion of Theorem \ref{theo1} applies.

\end{proof}

Given a contact metric structure $(M,\eta ,\xi ,\phi ,g)$, and a real constant $a>0$, then the $D_a$-homothetic deformation of $(M,\eta , \xi , \phi ,g )$ is another contact metric structure $(M,\overline{\eta}, \overline{\xi }, \overline{\phi}, \overline {g})$ where $\overline{\eta}=\frac{\eta}{a}$, $\overline{\xi}=\frac{\xi}{a}$, $\overline{\phi}=\phi$, $\overline{g}=ag+a(a-1)\eta\otimes\eta $.    Any  $D_a$-homothetic deformation of a generalized  Eta-Einstein K-contact structure is again Eta-Einstein. This is true even for strict generalized structures in dimension 3. To see this, one uses identity (\ref{12w}) to deforms the generalized Eta-Einstein structure as a $\kappa$-structure with constant $\kappa$, then use identity (\ref{eqric}) to see  that the deformed $\kappa$-structure is itself  a deformed generalized Eta-Einstein structure.   The Ricci operator transforms as follows under $D_a$-homothety \cite{SHA}: $$\overline{Q}X=\overline{\lambda}X+\overline{\gamma}\overline{\eta}(X)\overline{\xi}$$ where $$\overline{\lambda}=\frac{\lambda +2-2a}{a};~~~~~~~~\overline{\gamma}=2n-\overline{\lambda}$$

Eta-Einstein K-contact structures with $\lambda =-2$ are D-homothetically fixed. Indeed, any $D_a$-homothetic deformation of such an $\eta$-Einstein structure has the same $\overline{\lambda} =-2$. Any Eta-Einstein K-contact structure with $\lambda  >-2$ is D-homothetic to a K-contact, Einstein structure which usually have large isometry groups.  For those Eta-Einstein K-contact structures with $\lambda \le -2$, the isometry group tends to be small dimensional as evidenced through the following theorem.

\begin{Theorem} Let $(M,\eta ,\xi ,\phi ,g)$ be a closed, K-contact, generalized Eta-Einstein space whose structure is not D-homothetic to a K-contact, Einstein structure. Then the isometry group of $(M,\eta ,\xi ,\phi ,g)$ is one dimensional. In particular, its Reeb flow is almost regular.\end{Theorem}

\begin{proof}  Let $Z$ be a Killing vector field on the generalized  Eta-Einstein space $(M,\eta , \xi ,\phi , g)$. We will show that $Z$ is a (constant) multiple of $\xi$.  Since $\gamma>2n+2>0$, it follows from Theorem \ref{theo1} that $Z$ is an infinitesimal automorphism of the generalized Eta-Einstein structure. Now, if $\overline{g}=ag+a(a-1)\eta\otimes\eta$, $\overline{\eta}=a\eta$,  is a $D_a$-homothetic deformation, then $Z$ is still an infinitesimal automorphism of the deformed structure. The Ricci tensor $\overline{Ric}$ of the deformed metric satisfies $$\overline{Ric}(X,Y)=\overline{\lambda}\overline{g}(X,Y)+\overline{\gamma}\overline{\eta}(X)\overline{\eta}(Y)$$  where, \cite{SHA} $$\overline{\lambda}=\frac{\lambda +2-2a}{a},~\overline{\gamma}=2n-\frac{\lambda +2-2a}{a}.$$ Moreover, $\overline{\lambda}\le -2$ also, since the property ''$\lambda \le -2$" is preserved under D-homothetic deformations on generalized $\eta$-Einstein K-contact spaces. Computing $\overline{Ric}(Z,Z)$, we find:

$$\begin{array}{rcl}\overline{Ric}(Z,Z)&=&(\lambda +2-2a)g(Z,Z)+[2na^2+2a-(\lambda +2)]\eta (Z)^2\end{array} $$ which can be reorganized as: 
\begin{equation}\label{id1}\overline{Ric}(Z,Z)=(\lambda +2-2a)[g(Z,Z)-\eta (Z)^2]+2na^2\eta (Z)^2\end{equation}

From identity (\ref{id1}), we deduce that {\small  \begin{equation}\label{eqf1}\overline{Ric}(Z,Z)\le -2a[g(Z,Z)-\eta(Z)^2]+2na^2\eta(Z)^2\le -2a[g(Z,Z)-\eta(Z)^2-na\eta(Z)^2]\end{equation}}
If $\eta(Z)^2\ne g(Z,Z)$, that is, if $Z$ is not proportional to $\xi$, then  inequality (\ref{eqf1}) clearly implies that, for sufficiently small $a>0$:   \begin{equation}\label{inf1}\overline{Ric}(Z,Z)\le 0.\end{equation}

 On a closed Riemannian manifold, inequality (\ref{inf1}) would imply that the Killing vector field $Z$ is parallel, hence trivial  since there are no non-trivial parallel vector fields on any closed K-contact manifold \cite{RUK}.

 \end{proof}

\section{Generalized Jacobi $(\kappa ,\mu )$-spaces}
It was pointed out that examples of non-compact generalized $(\kappa ,\mu )$-spaces were provided in \cite{KOC}. However, no compact examples of those generalized $(\kappa ,\mu )$-spaces are      known. There is a weaker notion of $(\kappa ,\mu )$-structures, the so-called Jacobi $(\kappa ,\mu )$ structures \cite{GHS}. A contact metric structure $(M,\eta ,\xi ,J,,g)$ is called a generalized Jacobi $(\kappa ,\mu )$ if the Riemann curvature tensor satisfies
\begin{equation}R(X,\xi )\xi =\kappa (X-\eta(X)\xi)+\mu hX\end{equation} where $h=\frac{1}{2}L_\xi J$ and $\kappa$ and $\mu$ are smooth functions. Similar to the standard structures, for $\kappa <1$, the tangent bundle decomposes into 3 subbundles, the eigenbundles of $h$ corresponding to eigenfunctions $\lambda =\sqrt{1-\kappa }$, $-\lambda$ and $0$. We focus on the 3-dimensional case. Starting with a closed, flat, 3-dimensional  contact metric structure, let $E$ be a unit vector field such that $hE=E$. Then, the following identities hold (\cite{BKP}).
\begin{equation}\nabla E=0,~[\xi ,E]=2\phi E,~[\xi ,\phi E]=0,~[E,\phi E]=2\xi\label{17}\end{equation}
The 3-dimensional contact metric manifold fibers over the circle $S^1$ with $\xi$ and $\phi E$ tangent to the fibers. Therefore, there are plenty of functions invariant along the fibers, that is, functions $f$ such that $df(\xi )=0$ and $df(\phi E)=0$.

Any of these functions induces a contact metric deformation as follows:

Define a new metric $g^f$ by 
\begin{equation}g^f(E,E)=1+f+\frac{1}{2}f^2,~~g^f(\phi E,\phi E)=1-f+\frac{1}{2}f^2\end{equation}
\begin{equation} g^f(E,\phi E)=f^2, ~~g^f(\xi ,E)=0=g^f(\xi ,\phi E)\end{equation}
One defines transverse  almost complex structure $\phi^f$ as follows: For any tangent vectors  $X$, $Y$, \begin{equation}d\eta (X,Y)=2g^f(X,\phi^f Y)\end{equation}

A simple calculation shows immediately that \begin{equation}\phi^f\xi =0\end{equation}
\begin{equation}\phi^f E=-f^2E+(1+f+\frac{1}{2}f^2)\phi E\end{equation}
\begin{equation}\phi^f \phi E=-(1-f+\frac{1}{2}f^2)E+f^2\phi E\end{equation}
It is easily verified that $(M,\eta , \xi ,\phi^f ,g^f)$ is a contact metric structure with $h^f=\frac{1}{2}L_\xi \phi^f$ satisfying 
\begin{equation} h^f\phi E=-(1-f+\frac{1}{2}f^2)\phi E\label{23}\end{equation}
We prove the following:
\begin{pro} With the same notations as above, $(M,\eta , \xi ,\phi^f, g^f)$ is a generalized Jacobi $(\kappa ,\mu )$-structure with $\kappa =(f-\frac{1}{2}f^2)(2-f+\frac{1}{2}f^2)$ and $\mu =2(f-\frac{1}{2}f^2)$.
\end{pro}
\begin{proof} We will use the notation $\nabla$ instead of the cumbersome $\nabla^f$ for the Levi-Civita covariant derivative of $g^f$.  A direct calculation using identities (\ref{17}) and (\ref{23}) shows that the curvature tensor $R^f$  of $g^f$ satsfies:
$$\begin{array}{rcl} R^f(\phi E,\xi )\xi&=&-\nabla_\xi\nabla_{\phi E}\xi-\nabla_{[\phi E,\xi ]}\xi\\&=&\nabla_\xi (\phi^f \phi E+\phi^f h^f\phi E)\\&=&(f-\frac{1}{2}f^2)\nabla_\xi \phi^f \phi E\\&=&(f-\frac{1}{2}f^2)\phi^f\nabla_{\phi E}\xi\\&=&(f-\frac{1}{2}f^2)^2\phi E\\&=&(f-\frac{1}{2}f^2)(2-f+\frac{1}{2}f^2)\phi E-2(f-\frac{1}{2}f^2)(1-f+\frac{1}{2}f^2)\phi E\\&=&\kappa \phi E+\mu h^f\phi E
\end{array}$$
Similarly,
$$\begin{array}{rcl}R^f(\phi^f \phi E,\xi )\xi&=&-\nabla_\xi\nabla_{\phi^f \phi E}\xi-\nabla_{[\phi^f \phi E,\xi ]}\xi\\&=&-\nabla_\xi\nabla_{\phi^f \phi E}\xi -2(1-f+\frac{1}{2}f^2)\nabla_{\phi E}\xi\\&=&-\nabla_\xi (2-f+\frac{1}{2}f^2)\phi E+2(1-f+\frac{1}{2}f^2)(f-\frac{1}{2}f^2)\phi^f \phi E\\&=&(2-f+\frac{1}{2}f^2)(f-\frac{1}{2}f^2)\phi^f \phi E+2(1-f+\frac{1}{2}f^2)(f-\frac{1}{2}f^2)\phi^f \phi E\\&=&\kappa \phi^f\phi E+\mu (1-f+\frac{1}{2})\phi^f \phi E \\&=&\kappa\phi^f\phi E+\mu h^f\phi^f \phi E
\end{array}
$$

$\qed$
\end{proof}
\begin{Remark} The above construction does not seem to provide any example of generalized $(\kappa ,\mu )$-structures as we have hoped for. This failure is mainly due to the fact that the curvature tensor $R^f$ does not satisfy this condition \begin{equation} R^f(\phi E, \phi^f\phi E)\xi=0\label{cond}\end{equation}
 which is necessary on a  generalized $(\kappa ,\mu )$-structure. In our case at hand, condition (\ref{cond}) can be satisfied only if the function $f$ is  constant.
\end{Remark}

\end{document}